\documentclass[12pt,twoside,reqno]{amsart}

\usepackage{color}
\usepackage{graphicx}
\usepackage{amsmath}
\usepackage{amssymb}
\usepackage{amsfonts}
\usepackage{multicol}
\usepackage[ansinew]{inputenc}
\usepackage{amscd}
\usepackage[mathscr]{eucal}
\usepackage{hyperref}
\usepackage{xcolor}
\usepackage{easyReview}

\newcommand{\blob}{\rule[.2ex]{.8ex}{.8ex}}

\newcommand{\R}{\mathbb{R}}
\newcommand{\C}{\mathbb{C}}
\newcommand{\N}{\mathbb{N}}
\newcommand{\Z}{\mathbb{Z}}

\newcommand{\Su}{\mathbb{S}}
\newcommand{\SL}{{\rm SL}}

\newcommand{\GL}{\rm GL}

\newcommand{\Mat}{{\rm Mat}}

\newcommand{\abs}[1]{\bigl| #1 \bigr|} 
\newcommand{\norm}[1]{\lVert#1\rVert} 
\newcommand{\normtwo}[1]{
	{\left\vert\kern-0.25ex\left\vert\kern-0.25ex\left\vert #1
		\right\vert\kern-0.25ex\right\vert\kern-0.25ex\right\vert} }

\newcommand{\Proj}{\mathbb{P}(\R^ 2)}

\newcommand{\Pp}{\mathbb{P}}

\newcommand{\medir}{{\overline{v}}}
\newcommand{\ledir}{{\underline{v}}}

\newcommand\restr[2]{{
		\left.\kern-\nulldelimiterspace 
		#1 
		\vphantom{\big|} 
		\right|_{#2} 
}}

\theoremstyle{plain}
\newtheorem{theorem}{Theorem}[section]
\newtheorem{proposition}{Proposition}[section]
\newtheorem{corollary}[proposition]{Corollary}
\newtheorem{lemma}[proposition]{Lemma}

\theoremstyle{definition}
\newtheorem{definition}{Definition}[section]
\newtheorem{question}{Question}[section]

\theoremstyle{definition}
\newtheorem{remark}{Remark}[section]
\newtheorem{example}[theorem]{Example}
\numberwithin{equation}{section}

\newcommand{\Mcal}{\mathcal{M}}

\newcommand{\Prob}{\mathrm{Prob}}

\newcommand{\supp}{\mathrm{supp}}

\newcommand{\Tscr}{\mathscr{T}}
\newcommand{\Cscr}{\mathscr{C}}

\newcommand{\bld}[1]{\mathbf{{#1}}}

\title [Obstructions to the regularity of the Lyapunov exponents]{Obstructions to the regularity of the Lyapunov exponents for non-compact random Schr\"odinger cocycles}

\date{}

\begin{document}

\author[P. Duarte]{Pedro Duarte}
\address{Departamento de Matem\'atica and CEMS.UL\\
Faculdade de Ci\^encias\\
Universidade de Lisboa\\
Portugal
}
\email{pmduarte@fc.ul.pt}

\author[T. Graxinha]{Tom\'e  Graxinha}
\address{Departamento de Matem\'atica and CEMS.UL\\
	Faculdade de Ci\^encias\\
	Universidade de Lisboa\\
	Portugal}
\email{graxinhatome@hotmail.com}

\begin{abstract}
In this paper, we present a class  of random Schr\"odin\-ger cocycles  showing that, 
for random  cocycles with non-compact support, 
the presence of certain finite moment conditions is essential 
for establishing a specific modulus of continuity of the Lyapunov exponent. In particular, H\"older continuity of the Lyapunov exponent requires an exponential moment condition.
\end{abstract}

\maketitle


\section{Introduction}

A linear cocycle with values in $\SL_m(\R)$ over a measure preserving dynamical system $(X,\mathcal F,\mu,f)$ is a bundle map $F:X\times\R^m\to X\times\R^m$ of the form $F(\omega,v)=(f(\omega),A(\omega)v)$, where $A:X\to \SL_m(\R)$ is measurable; its $n$-th iterate is given by $F^n(\omega,v)=(f^n(\omega),A^{n}(\omega)v)$ with $A^{n}(\omega):=A(f^{n-1}(\omega))\cdots A(\omega)$. The (top) Lyapunov exponent is defined as
\[
L_1(F,\omega)\;:=\;\lim_{n\to\infty}\frac1n\log\|A^{n}(\omega)\|,
\]
whenever the limit exists. Furstenberg and Kesten proved (see ~\cite{FK60}) that $L_1$ exists and is finite  $\mu$-a.e.\ under the following first moment condition
$$\int\log \norm{A(\omega)} \, d\mu(\omega)<\infty.$$ Moreover, if the base dynamics is ergodic, $L_1(F,\cdot)$ is almost surely constant. Oseledets' theorem (see \cite{Ose68}) then yields the full spectrum  $L_1\geq \cdots \geq L_m$   and their respective multiplicities.
The second Lyapunov exponent is given by $$L_2(F;\omega):=\lim_{n\rightarrow+\infty}\frac{1}{n}\log s_2(A^n(\omega)),$$ where $s_2$ is the second singular value. A probability measure $\mu$ on $\SL_m(\R)$ defines a random cocycle $$F(\omega,v):=(\sigma\omega,A(\omega)v),$$ over the Bernoulli shift $\sigma\colon X\to X$ on $X=\SL_m(\R)^\Z$, endowed with the Bernoulli measure $\mu^\Z$, with the locally constant fiber action $A(\omega)=\omega_0$. The corresponding first and second Lyapunov exponents are denoted $L_1(\mu)$ and $L_2(\mu)$.
Thus, one can identify a random cocycle with a probability measure $\mu\in\Prob(\SL_m(\R))$.
When $supp(\mu)$ is compact, the finite moment condition is trivially satisfied. Furstenberg's positivity criterion (see \cite{Fur63}) implies that $L_1(\mu)>0$ whenever the semigroup generated by $\supp(\mu)$ is non-compact and strongly irreducible. Furstenberg and Kifer (see \cite{FK83}) established the generic continuity of the Lyapunov exponent, i.e. under irreducibility and a uniform first moment conditions. Also, under irreducibility, a spectral gap ($L_1>L_2$) and a uniform exponential moment, Le Page proved in~\cite{LP89} the H\"older continuity of $L_1$ for one-parameter  families of random cocyles.
Duarte and Graxinha (see \cite{DG25}) obtained the H\"older continuity of $L_1$ in general spaces of measures on $\Mat_m(\R)$ under the same hypothesis of Le Page, i.e.  finite exponential moment, quasi-irreducibility and a spectral gap.  Other results concerning the regularity of the Lyapunov exponent include the following:~\cite{DK-20} establish weak H\"older continuity under a spectral gap assumption and for finitely supported cocycles; 
	while \cite{Tall-Viana} prove pointwise H\"older continuity in the presence of a spectral gap, and pointwise log-H\"older continuity in general. 

 For compactly supported measures, general continuity, without ge\-ne\-ric assumptions, was established by Bocker-Neto and Viana for $\GL_2(\R)$-cocycles (see \cite{BV}) and, in the broader $\GL_m(\R)$ setting,  by Avila, Eskin and Viana (see \cite{AEV}).

Random cocycles typically admit a unique stationary measure. 
Under a generic irreducibility assumption together with an exponential moment condition, Guivarc’h~\cite{Guivarch90} proved that this stationary measure is H\"older regular. 
More recently, Gorodetski et al.~\cite{GKM22} extended this result to smooth random dynamical systems generated by a probability measure on a group of diffeomorphisms of a closed Riemannian manifold. 
Most recently, Monakov~\cite{Mon25} established log-H\"older regularity of the stationary measure for Lipschitz random dynamical systems under a linear (logarithmic) moment-type condition.

\bigskip 

There is a well-known connection between the spectral theory of Schr\"odinger operators and the Lyapunov exponents of linear cocycles.

Consider an invertible ergodic transformation $f:X\to X$ over a probability
space $(X,\mu)$. Given a bounded and measurable observable $v:X\to\mathbb{R}$,
let $v_n(\omega):=v(f^n \omega)$ for all $\omega\in X$ and $n\in\mathbb{Z}$.

\medskip
Denote by $\ell^2(\mathbb{Z})$ the Hilbert space of square summable sequences
of real numbers $(\psi_n)_{n\in\mathbb{Z}}$. The discrete ergodic Schr\"odinger
operator with potential $n\mapsto v_n(\omega)$ is the operator $H_\omega$ defined on
$\ell^2(\mathbb{Z}) \ni \psi=\{\psi_n\}_{n\in\mathbb{Z}}$ by
\begin{equation}
  (H_\omega\psi)_n := -(\psi_{n+1}+\psi_{n-1}) + v_n(\omega)\,\psi_n .
\end{equation}

Note that due to the ergodicity of the system, the spectral properties of the
family of operators $\{H_\omega:\omega\in X\}$ are $\mu$-a.s.\ independent of the phase
$\omega$.

\medskip
Given an energy parameter $E\in\mathbb{R}$, the Schr\"odinger (or eigenvalue)
equation $H(\omega)\psi=E\psi$ can be solved formally by means of the iterates of a
certain dynamical system. More precisely, consider the associated Schr\"odinger
cocycle $X\times\mathbb{R}^2\to X\times\mathbb{R}^2$,
$(\omega,v)\mapsto (f(\omega), A_E(\omega)v)$, where $A_E:X\to \mathrm{SL}_2(\mathbb{R})$ is
given by
\[
  A_E(\omega):=
  \begin{bmatrix}
    v(\omega)-E & -1\\[2pt]
    1      & 0
  \end{bmatrix}
  =
  \begin{bmatrix}
    v(\omega) & -1\\[2pt]
    1    & 0
  \end{bmatrix}
  \;+\;
  E\begin{bmatrix}
    -1 & 0\\[2pt]
     0 & 0
  \end{bmatrix}.
\]
Let $A_E^{n}$ denote the $n$-th iterate of the cocycle, that is,
\[
  A_E^{n}(\omega)=A_E(f^{\,n-1}\omega)\cdots A_E(f(\omega))A_E(\omega).
\]
Then the formal solution of the Schr\"odinger equation $H(\omega)\psi=E\psi$ is given by
\begin{equation}
  \begin{bmatrix}\psi_n\\ \psi_{n-1}\end{bmatrix}
  = A_E^{(n)}(\omega)\begin{bmatrix}\psi_0\\ \psi_{-1}\end{bmatrix}.
\end{equation}

\medskip
The top Lyapunov exponent of the Schr\"odinger cocycle is denoted by $L_1(A_E)$.

 By ergodicity, the spectra of these family of operators is almost surely constant  (see \cite{Pas80}).
Johnson's theorem (see \cite[Theorem 3.12]{David-survey}) establishes that the spectrum's complement corresponds to parameters where the Schr\"odinger cocycle is uniformly hyperbolic.

The integrated density of states (IDS) is a distribution function $N(E)$ that physically measures how many states
correspond to energies less than or equal to $E$. Mathematically, this corresponds to the asymptotic distribution of the eigenvalues of increasingly large Schrödinger matrices obtained by truncating the Schrödinger operator.
The Thouless formula (see~\cite[Theorem 3.16]{David-survey})
relates the Lyapunov exponent with the IDS
 $$L_1(\mu_E)=\int \log |E-E'| dN(E') , $$ 
 expressing it as the Hilbert transform of   $N(E)$. This formula was  initially employed by Craig and Simon (see \cite{CS83}) to prove the log-H\"older continuity of the IDS).  
 A threshold for the regularity preserved under the Hilbert transform was established by Goldstein and Schlag \cite[Lemma~10.3]{GS01}. For example, H\"older regularity lies above this threshold, whereas log-H\"older regularity falls below it.  More recently, Avila et al.\ \cite[Proposition~2.2 and Corollary~2.3]{ALSZ21} improved upon the result of Goldstein and Schlag, showing that log-H\"older moduli of continuity  are not preserved but are instead mapped into lower  log-H\"older type of  regularity within the same family.

   The regularity of the IDS for the Anderson model has been the subject of several studies: explicit estimates were obtained in ~\cite{HV17} for the case of small coupling constants and in ~\cite{DG13} for Fibonacci Hamiltonians. Higher-order regularities of the IDS under appropriate hypotheses were obtained by ~\cite{ST85} and ~\cite{CK86}.

	Moreover, in~\cite{GorK25} Gorodetski et al proved that the (fibered) rotation number of any cocycle (including Schr\"odinger) is log-H\"older, which via a simple connection to the IDS provides an alternative dynamical proof of the IDS's log-H\"older regularity.

In \cite{BCDFK} Bezerra et al established an abstract dynamical Thouless-type formula for affine families of $\GL_2(\R)$ cocycles. Here, the IDS admits a dynamical description as the fibered rotation number.  More precisely, if $K_n(\omega, E_0)$ denotes the number of full turns in $\mathbb{P}^1$ performed by the projective curve
\[
(-\infty, E_0]\ni E\ \longmapsto\ A^{n}_{E}(\omega)\,\hat v,
\]
for a typical $\omega$ and any $\hat v\in\mathbb{P}^1$, 
\[
N(E_0)\;=\;\lim_{n\to\infty}\frac{K_n(\omega, E_0)}{n},
\]
and this rotation number agrees with the IDS (see~\cite[Lemma 3.14]{BCDFK} or~\cite[Section 3]{GorK25}).

Let $\mu$ be a probability measure on the real line.  
	The two-sided Bernoulli shift $\sigma:\R^\Z \to \R^\Z$, endowed with the product measure $\mu^\Z$, is a classical example of an ergodic and mixing measure-preserving dynamical system.  
	Consider the locally constant function $v:\R^\Z \to \R$ defined by $v(\omega) = \omega_0$,  
	which generates a random i.i.d.\ potential via
	\[
	v_n(\omega) := v(\sigma^n \omega) = \omega_n.
	\]
	The measure $\mu$ determines the random Schr\"odinger cocycle 
	 $\mu_E=(A_E)_\ast \mu^\Z$,   also  characterized as
	\[
 \mu_E = \int_{-\infty}^\infty 
	\delta_{\begin{bmatrix}
			x -E & -1 \\ 
			1 & 0
	\end{bmatrix}} \, d\mu(x). 
	\]

	Every random Schr\"odinger cocycle $\mu_E$ is strongly irreducible and non-compact (see \cite[Subsection~4.3]{David-survey}).  
	By Furstenberg's criterion, this implies that the top Lyapunov exponent is positive.  
	Moreover, since Schr\"odinger cocycles take values in $\SL_2(\R)$, the Lyapunov spectrum exhibits a gap:
	\[
	L_1(\mu_E) > 0 > -\,L_1(\mu_E) = L_2(\mu_E).
	\]

\bigskip

	In~\cite{BCDFK} (see also~\cite{BD23}), in the setting of random strongly irreducible cocycles generated by finitely supported measures, the authors establish an upper bound for the H\"older regularity of the Lyapunov exponent. This bound coincides with the ratio between Kolmogorov entropy of the underlying shift and the Lyapunov exponent. Typically, this ratio equals twice the dimension of the stationary measure. 
	   
	Other negative results regarding the regularity of the IDS were obtained in ~\cite{GK11},	which proves the optimality of the log-H\"older continuity of the IDS for certain examples of dynamically defined Schr\"odinger operators and in ~\cite{M19} with a class of quasi-periodic Schr\"odinger operators where the IDS is not H\"older continuous for a sufficiently large coupling constant.

	Recall that a probability measure $\eta$ on $\Pp^1$ is said to be $\alpha$-H\"older if
	\[
	\eta(B(x,r)) \leq C\, r^\alpha \quad \forall\, x \in \Pp^1 \ \forall\, r>0,
	\]
	which in particular implies that $\dim(\eta) \geq \alpha$.

	An interesting new  connection between the regularity of a cocycle's stationary measure and that of its rotation number was recently established in~\cite{DGK26}. Via the standard relation between the rotation number of a Schrödinger cocycle and the IDS, this provides an alternative way of  studying the regularity of the IDS.

	These observations suggest that, at least in this context, the following three quantities are intimately related:
	\begin{itemize}
		\item the regularity of the Lyapunov exponent around $\mu$;
		\item the regularity of the stationary measures of $\mu$;
		\item the dimension of the stationary measures of $\mu$.
	\end{itemize}

\bigskip

In this work we construct a random unbounded Schr\"odinger cocycle with locally uniformly bounded \emph{sub-exponential moments}
\[
\sup_{|E|\le m}\int \exp\!\big((\log\|g\|)^{1/3}\big)\,d\mu_E(g)\;<\;\infty
\quad\text{for every }m>0,
\]
but with \emph{infinite exponential moments}
\[
\int \|g\|^\alpha\,d\mu_E(g)\;=\;\infty
\quad\text{for every }\alpha>0 \; \text{ and } \; E,
\]
such that the Lyapunov exponent $E\mapsto L_1(\mu_E)$ is not locally $\alpha$-H\"older continuous for any $\alpha>0$. In particular, this shows that all the hypothesis in \cite{LP89} and \cite{DG25}, except for the exponential moment condition, are satisfied, thereby demonstrating the sharpness of these results.

More generally, we set up a dictionary between \emph{moment profiles},
see definitions~\ref{def moment profile} and~\ref{def moment satisfied}, and \emph{moduli of continuity} such that when a given (locally uniform) moment condition fails, the corresponding modulus of continuity for $E\mapsto L_1(\mu_E)$ cannot hold (see Theorem~\ref{main}).

This raises a natural question: 
Does there exist a $1$-$1$ correspondence $\varphi\leftrightarrow\omega$ between moment profiles and moduli of continuity such that, whenever a random Schr\"odinger cocycle $(\mu_E)_E$ satisfies a moment profile $\varphi_0$ locally uniformly in $E$ then the map $E\mapsto L_1(\mu_E)$ satisfies the associated local modulus of continuity $\omega_0$? 
An analogous question can be posed for the dependence of $L_1$ on the generating law $\mu\in\mathrm{Prob}(\mathrm{SL}_m(\mathbb{R}))$ with respect to the Wasserstein distance, under the usual irreducibility and spectral gap assumptions.

Similar questions may be posed regarding the relation between moment conditions and the regularity of the corresponding stationary measures, to which the works~\cite{GKM22,Mon25} already provide partial affirmative answers.

A positive answer to these questions would significantly clarify the picture on the quantitative regularity of Lyapunov exponents in the non-compact settings. We note that, for compactly supported random $\mathrm{SL}_m(\mathbb{R})$ cocycles, generic H\"older dependence on $\mu$ with respect to the Wasserstein distance is known.

\section{Main Results and questions}

\begin{definition}
A function \(\omega \colon [0, 1) \to [0, 
+\infty)\) is called a \textit{modulus of continuity} (MOC) provided it is:  
(i) continuous,  
(ii) strictly-increasing and 
(iii) \(\omega(0) = 0\).  
\end{definition}

Let \((X, d)\) be a metric space.
\begin{definition}
	A function \(f \colon X \to \mathbb{R}\) is said to have \textit{local modulus of continuity} \(\omega\) if for every $a\in X$, there exist positive constants $r>0$ and $C<\infty$ such that for all $x,y\in X$ with
	$d(x,a)<r$ and $d(y,a)<r$,
\begin{equation*}
    |f(x) - f(y)| \leq C\, \omega(d(x,y)) .
\end{equation*}
\end{definition}

Common examples of Moduli of continuity are the following:

\noindent
\blob\; 
\textbf{H\"older continuity}.  
A function \(f \colon X \to \mathbb{R}\) is $\alpha$-\textit{H\"older continuous} if it has modulus of continuity   
\begin{equation}
    \omega(r) =r^\alpha = \exp\left(-\alpha \log \tfrac{1}{r}\right),
\end{equation}
where $0<\alpha\leq 1$.
The case \(\alpha = 1\) corresponds to \textit{Lipschitz continuity}.

\noindent
\blob\; \textbf{Weak-H\"older continuity}. 
A function \(f\) is $(\alpha, \theta)$-\textit{weak-H\"older continuous} if it has modulus of continuity   
\begin{equation}
    \omega(r) = \exp\left(- \alpha \left( \log \tfrac{1}{r}\right)^\theta\right),
\end{equation} 
for some $\alpha>0$ and  $0<\theta\leq 1$.  
When \(\theta = 1\), this coincides with H\"older continuity.

\noindent
\blob\; \textbf{Log-H\"older continuity}. 
A function \(f\) is $\gamma$-\textit{$\log$-H\"older continuous} if it has modulus of continuity
\begin{equation}
	\label{log Holder MOC}
    \omega(r) =  \left(\log \tfrac{1}{r}\right)^{-\gamma},
\end{equation}
where $\gamma>0$.

Moduli of continuity are partially ordered by the following relation: we say that $\omega'$ is finer than $\omega$, or that $\omega'$ implies $\omega$, and write  $\omega'\leq\omega$,  if there exists $C<\infty$ and $r_0>0$ such that 
$\omega'(r)\leq C \omega(r)$ for all $0<r\leq r_0$.
The previous classes are  hierarchies of MOC
each one ordered by its own parameter, larger parameters corresponding to finer  MOC.
The three classes are related as follows:
\begin{equation}
    \text{H\"older} \Rightarrow \text{weak-H\"older} \Rightarrow \text{log-H\"older}.
\end{equation}

\bigskip

\begin{definition}
	\label{def moment profile}
	A function \(\varphi \colon (1, +\infty) \to (0, 
	+\infty)\) is called a \textit{moment profile} provided it is:  
	(i) continuous,  
	(ii) strictly-increasing and satisfies\, 
	(iii) $\lim_{r\to \infty} \varphi(r)=\infty$.  
\end{definition}

\begin{definition}
	\label{def moment satisfied}
	Given a moment profile $\varphi$, we say that a measure $\mu\in \Prob(\SL_m(\R))$ \textit{has finite $\varphi$-moment}  if 
	\begin{equation*}
		\int \varphi(\log \norm{g})\, d\mu(g)<\infty . 
	\end{equation*}
\end{definition}

Common examples of moment profiles are the following:

\noindent
\blob\; 
\textbf{Exponential moment}.
We say that $\mu\in \Prob(\SL_m(\R))$ has \textit{finite exponential moment} if it has finite  moment profile   
\begin{equation}
	\varphi(r) = \exp\left(\alpha  r \right) \;  \text{ with } \quad \alpha>0 .
\end{equation}

\noindent
\blob\; \textbf{Sub exponential moment}.
We say that $\mu\in \Prob(\SL_m(\R))$ has \textit{finite sub-exponential moment} if it has finite  moment profile   
\begin{equation}
	\varphi(r) = \exp\left( r^\theta \right) \;  \text{ with } \; 0<\theta\leq 1.
\end{equation}

\noindent
\blob\; \textbf{Polynomial moment}. 
We say that $\mu\in \Prob(\SL_m(\R))$ has \textit{finite polynomial  moment} if it has finite  moment profile   
\begin{equation}
	\varphi(r) = r^\gamma\;  \text{ with } \quad \gamma>0 .
\end{equation}

Moment profiles are partially ordered by the following relation: we say that $\varphi'$ is stronger than $\varphi$, or that $\varphi'$ implies $\varphi$, and write  $\varphi'\geq\varphi$,  if there exists $C<\infty$ and $r_0> 1$ such that 
$\varphi(r)\leq C\, \varphi'(r)$ for all $r\geq r_0$.
The previous classes are  hierarchies of moment profiles
each one ordered by its own parameter, larger parameters corresponding to stronger  moment profiles.
The three classes are related as follows:
\begin{equation}
	\text{Exponential} \Rightarrow \text{Sub-exponential} \Rightarrow \text{Polynomial}.
\end{equation}

\bigskip

Given $\beta>0$ we define a bijective  transformation $\Tscr_\beta$
between the spaces of moment profiles and  of moduli of continuity, $\Tscr_\beta (\varphi) := \omega$,
\begin{equation}
\omega(r)= \frac{1}{\varphi(\log \frac{1}{r})^\beta} ,
\end{equation}  
whose inverse transformation $\varphi=\Tscr_\beta^{-1}(\omega)$ is given by
\begin{equation}
	\varphi(r)= \frac{1}{\omega(e^{-r})^{1/\beta}} .
\end{equation}

%
%
%
%

These maps will be used as dictionaries between finite moment conditions and moduli of continuity for the Lyapunov exponent.

\begin{lemma}
The bijection $\Tscr_\beta$ is order reversing  (stronger moment profiles correspond to finer MOC) and maps:
\begin{itemize}
	\item  $\alpha$-exponential moment profiles to  $\beta\alpha$-H\"older MOC;
	
	\item $\theta$-sub exponential moment profile to  $(\beta,\theta)$-weak H\"older MOC;
	
	\item  $\gamma$-polynomial moment profile  to  $\beta\gamma$-log H\"older MOC. 
\end{itemize}
\end{lemma}

Any matrix of the form 
	\begin{equation}
		\label{Sch mat}
		S_v := \begin{bmatrix}
			v & -1 \\ 1 & 0 
		\end{bmatrix}  
	\end{equation}
will be referred to as a  \textit{Schr\"odinger} matrix.

Our main results are the following:

\begin{theorem}
	\label{main}
Consider two moment profiles $\varphi, \psi$ such that 
$$r\leq \psi(r)\leq \varphi(r), \quad \forall\, r > 1 $$
  and  
the family of measures, $\mu_t \in \Prob(\SL_2(\mathbb{R}))$ 
\[
\mu_t = \sum_{n=1}^{\infty} \left[ \frac{p_n}{2} \delta_{A_{v_n,t}} + \frac{p_n}{2} \delta_{A_{-v_n,t}} \right], \quad t \in \mathbb{C}
\]
where

\begin{enumerate}
	\item  $\mu_t$ determines a Random Schr\"odinger cocycle with matrices 
	\[
	A_{v_n,t} = S_{ v_n - t} \quad \text{ and } \quad A_{-v_n,t} = S_{ -v_n - t} ; \]
	
	\item $\sum_{n\geq 1} p_n=1$ \, and \ $0<\limsup_{n\to \infty} \frac{p_{n-1}}{p_n}\leq 1$;
	
	\item  for every $n\geq 1$, $v_n > 0$, \ $v_{n+1}-v_n >4$   and  
	
	  \quad $\displaystyle \lim_{n\to \infty}  v_n =\lim_{n\to \infty}  v_n -v_{n-1}= \infty;$
	\item The measures $\mu_t$ have locally uniformly bounded $\psi$-moments.
	\item $	\displaystyle \lim_{n\to \infty} p_n\varphi\big(\log v_n \big)=\infty$, which implies that the measures $\mu_t$ do not have finite $\varphi$-moments.
\end{enumerate}
 Under these assumptions, the IDS function $\R\ni E \mapsto N(E)$, i.e., the fibered rotation number of the Schr\"odinger cocycle $A_E^n(\omega)$,
 can not have $\omega=\Tscr_2(\varphi)$    as a local  MOC. 

\end{theorem}

\begin{remark}
In the previous theorem, the IDS is considered the unique distribution function $\rho$ for which the Thouless formula holds; see Proposition~\ref{prop existence of ids}. The IDS is generally defined as the asymptotic distribution function of the eigenvalues of truncations of the Schrödinger operator on intervals of the one-dimensional lattice as the size of these intervals tends to infinity.
It is not difficult to verify that these two definitions agree. See, for example, Remark 6.8 of~\cite{BL} or the last sentence of~\cite{CS83}.
\end{remark}

\begin{corollary}
	\label{maincoro}
	
Under the hypothesis  of Theorem~\ref{main},
\begin{enumerate}
	\item[(a)] If $\varphi(r)$ is an exponential or sub-exponential moment profile then the Lyapunov exponent function $\R\ni E\mapsto L_1(\mu_E)$ does not have $\omega=\Tscr_2(\varphi)$    as a local  MOC. 
	
	\item[(b)]  If $\varphi(r)=r^p$ for some $p>1$,  then the function $\R\ni E\mapsto L_1(\mu_E)$  is not locally $(2p-1)$-log H\"older continuous.
	
\end{enumerate}
\end{corollary}

\bigskip

The previous Schr\"odinger cocycle is associated with the
unbounded  $1$-dimensional discrete Schr\"odinger operator $H_\omega$ on $\ell^2(\Z)$ defined by
\begin{equation}
	\label{SE}
 (H_\omega\zeta)_n:=-(\zeta_{n+1}+\zeta_{n-1})\, +  v(\sigma^n\omega)\, \zeta_n , 
\end{equation}
where $v(\omega):=\omega_0$ \ and   the sequence $\omega=(\omega_n)_{n\in \Z}$ (regarded as a process) is \ i.i.d. with respect to the Bernoulli measure $\nu^\Z$, where  
\begin{equation}
	\label{PotMeas}
\nu = \sum_{n=1}^{\infty} \left[ \frac{p_n}{2} \delta_{v_n} + \frac{p_n}{2} \delta_{-v_n} \right] .
\end{equation}
In other words \ 
\,  $ \Pp[ \, \omega_n= -v_j \, ]= \frac{ p_j }{2}=  \Pp[ \, \omega_n=v_j \, ] $ \ for all $j\geq 1$ 
and these events are independent for different sites $n\in\Z$.

\begin{corollary}\label{corollary}
The finite exponential moment hypothesis is essential for the H\"older regularity of the Lyappunov exponent in \cite{LP89} and \cite{DG25}.
\end{corollary}

Given a positive  $C<\infty$ and a moment profile $\varphi$ consider the space $\Mcal^\varphi_C$
of probability measures $\mu\in \Prob(\SL_m(\R))$
such that
$$\int \varphi(\log \norm{g})   \, d\mu(g) \leq C   .  $$

\bigskip

We can now formalize the main question stated in the introduction.

\begin{question}
Is there a constant $\beta > 0$ such that for any moment profile 
$\varphi(r) \geq r$, and for any quasi-irreducible 
$\mu \in \mathcal{M}_C^\varphi$ with $L_1(\mu) > L_2(\mu)$, 
the Lyapunov exponent
\[
\mathcal{M}_C^\varphi \ni \mu \;\longmapsto\; L_1(\mu)
\]
admits a local modulus of continuity 
$\omega = \mathcal{T}_\beta(\varphi)$ around $\mu$, 
with respect to the Wasserstein distance on $\mathcal{M}_C^\varphi$?
\end{question}

\begin{question}
\label{second question}
In Theorem~\ref{main}, does  the same lack of the local MOC $\omega$ \ hold  in the neighborhood of each energy value in the spectrum ?
\end{question}

\section{Proofs of the main results}

This section contains the proofs of Theorem~\ref{main} and its corollaries.

Because $\mu_t$ generates a random (non-constant)  Schr\"odinger cocycle, by~\cite[Subsection 4.3]{David-survey}, $\mu_t$ is non-compact, strongly irreducible and $L_1(\mu_t)>0$ for all $t\in\R$.

\begin{proposition}\label{p1}The family of measures $\{\mu_t : t\in\mathbb{C}\}$ has locally uniform finite first  moment.
\end{proposition}

\begin{proof}
	For each $t\in\mathbb{C}$
	\begin{align*}
		\int \log\norm{g} \ d\mu_t(g) \leq \int \psi(\log\norm{g})\, d\mu_t(g)  <+\infty.
	\end{align*}
	And since $\mathrm{supp}(\mu_t)\subset \SL_2(\C)$, we also have  $$\int \log\norm{g^{-1}} \, d\mu_t(g)= \int \log\norm{g} \, d\mu_t(g)<+\infty.$$
 Notice that by the singular value  decomposition, $\norm{g}=\norm{g^{-1}}$ is the first singular value of every $g\in\SL_2(\C)$. 
\end{proof}

The previous bounds and Kingman's Subadditive Ergodic Theorem imply that the Lyapunov exponent exists
$$L_1(\mu_t)=\lim_{n\rightarrow +\infty}\frac{1}{n}\int\log\norm{g} \ d\mu_t^n(g),$$
where for all $t\in \C$,  $\mu_t^n$ denotes the convolution $n$-th power of $\mu_t$.
Next proposition states the continuity and subharmonicity of the Lyapunov exponent as a function on the complex plane.
Let 
\[
I_n := [v_n-2, v_n+2], 
\qquad 
I_{-n} := [-v_n-2, -v_n+2].
\]

\begin{proposition}
	\label{prop CSHHn1}
	The function  $\C\ni t\mapsto L_1(\mu_t)$ 
	is 
	\begin{enumerate}
		\item continuous on $\C$;
		\item subharmonic on $\C$;
		\item harmonic on $\C\setminus \Sigma$, where \
		$ \Sigma:= \bigcup_{n=1}^\infty \left(\, I_n \cup I_{-n} \, \right)$.
	\end{enumerate}
\end{proposition}

\begin{remark}
 $\Sigma$  is the spectrum of the unbounded Schr\"odinger operator $H_\omega$ in~\eqref{SE} 
  for $\nu^\N$-almost every $\omega\in\R^\Z$,  where  $\nu$ was introduced in~\eqref{PotMeas}.  Compare with~\cite[Theorem 5.2.2]{DF-book}.
\end{remark}

\begin{proof}
By~\cite[Proposition 4.1]{FK83} and Proposition~\ref{p1}, the function $L_1(\mu_t)$ is continuous in $t$.

Item (2) is a well known fact that follows from ~\cite{CFKS87}[Section 9.4].

Item (3) holds because the random cocycle determined by $\mu_t$ is uniformly hyperbolic for all $t\in \C\setminus\Sigma$ and, by a classical result of D. Ruelle~\cite{Rue79a},
the function $t\mapsto L_1(\mu_t)$ is analytic and harmonic
for $t\in \C\setminus \Sigma$.

Notice that when $\mathrm{Im}(v)\neq 0$ the matrices $S_v$ strictly contract one of the hemispheres determined by $\R\Pp^1$ in $\C\Pp^1$,  and hence are hyperbolic (see the proof of ~\cite{Ycz04}[Lemma 2]).
Additionally, for $t\in \R\setminus \Sigma$, the matrices in $\supp(\mu_t)$ have the form $S_v$ with $v\in \R\setminus [-\lambda, \lambda]$ for some constant $\lambda>2$.  Since all these matrices strictly  contract a $45^\circ$ cone along the $x$-axis, the semigroup generated by $\supp(\mu_t)$ is hyperbolic. 
\end{proof}

We have seen that
$\mu_t$ is uniformly hyperbolic while $L_1(\mu_t)$ is analytic and harmonic for $t \in \C\setminus \Sigma$.
The same holds for the following truncated measure.

Given $N \in \mathbb{N}$, we consider the normalized truncated measure
\[
\mu_{N,t} = \left( \sum_{n=1}^{N} p_n \right)^{-1} \left( \sum_{n=1}^{N} \left[ \frac{p_n}{2} \delta_{A_{v_n, t}} + \frac{p_n}{2} \delta_{A_{-v_n,t}} \right] \right) \in \Prob(\SL_2(\R))
\]
and the Lyapunov exponent of the random cocycle it determines 
$$L_1(\mu_{N,t}):=\lim_{n\rightarrow\infty}\frac{1}{n}\int\log\norm{g} \ d\mu^n_{N,t}(g).$$

 This random linear cocycle is the Schrödinger cocycle associated with the Anderson operator  dynamically generated by the potential $v:\R^\Z\to\R$,
	$v(\omega):=\omega_0$, sampled along the orbits of the shift  $\sigma: \R^\Z\to\R^\Z$ endowed with the truncated Bernoulli measure $\nu_N^{\Z}$, where 
\begin{equation}
	\label{eq nuN}
	\nu_N = \Bigl( \sum_{j=1}^N p_j \Bigr)^{-1} 
	\sum_{j=1}^N \frac{p_j}{2}\, (\delta_{v_j} + \delta_{-v_j}).
\end{equation}
\begin{proposition}
	The Lyapunov exponents $t\mapsto L_1(\mu_{N,t})$	are
	\begin{enumerate}
		\item continuous on $\C$;
		\item subharmonic on $\C$;
		\item harmonic on $\mathbb{C} \setminus \Sigma_N$, where
		$ \Sigma_N:= \bigcup_{i=1}^N \left(\,  I_{-i} \cup I_i \, \right) $;
		\item $L_1(\mu_{t}) = \lim_{N\to \infty} L_1(\mu_{N,t})$,  for every $t\in \C$. Moreover, the 
		 convergence holds  uniformly over compact subsets $K\Subset \mathbb{C} \setminus \Sigma$.
	\end{enumerate}
\end{proposition}

\begin{proof}
Items (1)-(3) are well known in the theory of bounded Anderson Schr\"odinger operators. The first part of (4) is a consequence of~\cite[Theorem B]{FK83}. For the second part we use the mean value formula.
\end{proof}

Let $\rho_N:\R\to \R$ be the integrated density of states (IDS) of the operator $H_{N,\omega}$  determined by the truncated  Bernoulli measure $\nu_N^\Z$, which by~\cite{BCDFK}  or~\cite{GorK25} is also the fibered rotation number of the family of random cocycles $\mu_{N,t}$.
By the classical Thouless formula 
\begin{equation}
	\label{Thouless rhoN}
	L_1(\mu_{N,t} ) =  \int \log {|t - s|} \, d\rho_{N} (s). 
\end{equation} 
Next proposition serves to prove two statements: first the uniform convergence of $\rho_N$ to $\rho$, when $N\to\infty$, on compact intervals of the real line (see item (5))
and second  that the Thouless formula~\eqref{Thouless rhoN} holds for $\rho$ as well (see item (6)).
These are the key statements in the proof of Theorem~\ref{main}. From (4), we only need  the uniform moduli of continuity  of the truncated measures $\rho_N$.
Item (6) is also  a means of proving (5). 

\begin{proposition}\label{prop existence of ids}
	There exists $\rho:\R\to \R$ such that:
\begin{enumerate}
	\item $\rho$ is continuous;
	\item $\rho$ is non-decreasing;
	\item $\lim_{t\to -\infty} \rho(t)=0$,\;
	$\lim_{t\to +\infty} \rho(t)=1$;

\item  given $n\geq 1$, there exist constants $C_n<\infty$  such that
for all  $N\geq n$ and  all $t,s\in [-v_n-2, v_n+2]$ with $|t-s|\leq 1$, $$\abs{\rho_N(t)-\rho_N(s)}, \ \abs{\rho(t)-\rho(s)}  \leq \frac{C_n}{\log\frac{1}{ |t-s|}}.$$

\item $	\displaystyle \rho(t)=\lim_{N\to \infty} \rho_N(t)$, for all $t\in \R$, with uniform convergence on compact sets;

\item for all $t\in \C$.
$$ 	L_1(\mu_{t} ) =  \int \log {|t - s|} \, d\rho(s); $$

\end{enumerate}

\end{proposition}

\begin{proof}
We believe that the content of this proposition for Anderson models with unbounded potentials  satisfying a finite logarithmic moment hypothesis should be clear to many experts, but we were unable to find complete references for the proofs of all items.

Items (1)–(3) follow from \cite{CFKS87}[Section 9.4].

Item (4) follows from \cite{GorK25}[Theorem 3.1], where $\rho$ corresponds to both the fibered rotation number (of the Schr\"odinger cocycle) and the IDS (of the Anderson model).

Item (5) appears to follow from ~\cite[Theorem 1]{S21}, but a more detailed analysis shows that it does not. The problem is that for the Kantorovich-Rubinstein (Wasserstein) metric, a simple calculation gives that $\mathrm{d_{KR}}(\nu, \nu_N)=\infty$ and hence we can not conclude from this theorem that $\rho_N$ converges to $\rho$.

Item (6), for unbounded potentials with a finite logarithmic moment   is mentioned in several articles, for example, in the last sentence of ~\cite{CS83}, but we were unable to find a reference with a complete proof of this statement.

Hence, for the sake of completeness,  we decided to include the
intertwined proofs of items (5) and (6).

The strategy of the proof is first to show that (5) holds for a subsequence $\{N_j\}_j$ with a certain limit   $\rho$. Then we prove that  (6) (Thouless formula)   holds for this measure $\rho$. This implies that the sublimit $\rho$ is unique, it is the distributional Laplacian of the subharmonic function $\C\ni t\mapsto L_1(\mu_t)$. Finally the convergence in item (5) follows from the uniqueness of the sublimits.

By item (4), given any compact interval $I \subset \mathbb{R}$, the family $\{\rho_N\}_{N \in \mathbb{N}}$ is equicontinuous on $I$.  
By the Arzel\`{a}--Ascoli theorem, since $\{\rho_N\}_{N \in \mathbb{N}}$ is uniformly bounded with values in $[0,1]$, there exists a function $\rho_I: I \to \mathbb{R}$ and a subsequence $\{N_j\}_j$ such that
\[
\lim_{j \to \infty} \rho_{N_j}\vert_I = \rho_I
\quad \text{uniformly on   } I  .
\]
Taking a nested sequence of intervals   exhausting $\R$, by a diagonal 
argument there exists a non-decreasing function $\rho: \R \to \R$ and a (diagonal) subsequence $\{N_j\}_j$
such that
\[
\lim_{j \to \infty} \rho_{N_j} = \rho
\quad \text{uniformly on every compact interval } I \subset \R .
\]
This implies that
for any $t\in\C\setminus \R$ and $M > 0$,
\[
\lim_{N \to \infty} \int_{-M}^M \log|t-s| \, d\rho_N(s)
= \int_{-M}^M \log|t-s| \, d\rho(s).
\]
From this we get 
\begin{align}
	L_1(\mu_t)&= \lim_{j\to\infty} L_1(\mu_{N_j,t}) \nonumber\\
	&=\lim_{j\to\infty} \int \log|t-s|\, d\rho_{N_j}(s) \nonumber\\
	&=\int \log|t-s|\, d\rho(s)  \label{ThoulessFormula}
\end{align}
where in the last step we need to use the following
\begin{lemma}
	\label{lemma IDS remainder}
	For any $t\in\C\setminus \R$ and $M > 0$,
	$$ \lim_{M\to \infty}  \int_{\R\setminus [-M,M]} \log |t-s|\, d\rho_N(s) =0 $$ 
	uniformly in $N$. The same holds for $\rho=\lim_{j\to \infty} \rho_{N_j}$ instead of $\rho_N$.
\end{lemma}

\begin{proof}
	The functions $\rho_{N_j}$ are constant on 
	$\R \setminus \Sigma$.  Hence, taking the limit, the same holds for $\rho$. 
	Moreover, we claim that for $1 \leq |n| \leq N$,
	\begin{equation}
		\label{eq drhoN(In)}
		d\rho_N(I_n) = \Bigl( \sum_{j=1}^N p_j \Bigr)^{-1} \frac{p_n}{2}, 
	\end{equation}
	which passing to the limit as $N \to \infty$ yields
	\[
	d\rho(I_n) = \frac{p_n}{2}.
	\]
	
	Let us now prove~\eqref{eq drhoN(In)}.
	As explained in~\cite{BCDFK}, one may write
	\[
	d\rho_N(I_n) = \lim_{m \to \infty} 
	\frac{1}{\pi\, m}\, \ell_{I_n}\bigl(A_t^m(\omega) \hat v \bigr),
	\]
	where $\hat v \in \Pp^1$ is any projective point (for instance $\hat v=(1:0)$),  
	\[
	A_t^m(\omega) 
	= S_{ \omega_m - t}\, 
	\cdots  \,
	S_{\omega_1 - t}  \,
	S_{\omega_0 - t},
	\] 
	with the notation~\eqref{Sch mat}, $\omega=(\omega_j)_{j \geq 1}$ is a typical sequence for the measure $\nu_N^{\Z}$  in~\eqref{eq nuN}, and  $\ell_I(A_t \hat v)$ denotes the length of the projective curve $I \ni t \mapsto A_t \hat v \in \Pp^1$  (see Definition~\ref{def length} below).  This length divided by $\pi$ counts the number of full turns of the previous curve around $\Pp^1$.

	Now fix $n \geq 1$.  
	If $\omega_j \neq v_n$, then the matrix $A_t(\omega_j) := S_{ \omega_j - t }$ 
	remains hyperbolic   with large trace as $t$ ranges in $I_n$.  
	Thus $A_t(\omega_j)$ gives no full turn around $\Pp^1$  when $t$ varies in $I_n$.  
	On the other hand, when $\omega_j = v_n$ the trace of $A_t(v_n)$ varies from $-2$ to $2$ as $t$ ranges over $I_n$, and in this case $A_t(v_n)$ produces exactly one full turn on $\Pp^1$.  
	By~\cite[Proposition~2.18]{BCDFK}, this implies
	\[
	d\rho_N(I_n) \;\geq\; 
	\lim_{m \to \infty} \frac{1}{m} 
	\#\{\, 0 \leq j \leq m-1 : \omega_j = v_n \,\} 
	= \Bigl(\sum_{j=1}^N p_j \Bigr)^{-1} \frac{p_n}{2}  .
	\]
	By the hypothesis (3) of Theorem~\ref{main}, the intervals $I_n$ are pairwise disjoint.  Hence, since $d\rho_N$ is a probability measure,
	\[
	d\rho_N(I_n) = \Bigl( \sum_{j=1}^N p_j \Bigr)^{-1} \frac{p_n}{2} ,
	\]
	which proves the claim~\eqref{eq drhoN(In)}.

	Next, fix $t \in \C\setminus \R$ and let $M-2 > 2|t|$.  
	If $v_j \geq M-2$, then
	\[
	\max_{s\in I_{\pm j}}  |t-s| \leq | |t| - v_j - 2| = |v_j + 2 - |t||.
	\]
	Therefore
	\begin{align*}
		\left| \int_{\R \setminus [-M,M]} \log|t-s|\, d\rho(s)\right| 
		&\leq 2 \sum_{v_j \geq M-2} d\rho(I_j)\, \log|v_j+2-|t| | \\
		&\leq \sum_{v_j \geq M-2} p_j \left( \log|v_j| + \log\Bigl|1 + \tfrac{2-|t|}{v_j}\Bigr| \right) \\
		&\leq \sum_{v_j \geq M-2} p_j \left( \log|v_j| + \tfrac{2-|t|}{M-2} \right),
	\end{align*}
	which tends to $0$ as $M \to \infty$.  
	The same bounds apply to $\rho_N$, so the convergence is uniform in $N$.  
\end{proof}
		
By means of~\cite[Theorem 1.1]{CS83b}, Thouless's formula~\eqref{ThoulessFormula} is valid for all $t\in\C$.
This implies that the Lebesgue-Stieltjes measure $d\rho (t)$ is uniquely determined as the Laplacian (in the sense of distributions) of the subharmonic function
	$u(t)=L_1(\mu_t)$. 

Consequently, for every compact interval $I\subset \R$, every convergent sublimit of the sequence $\{\rho_N\vert_I\}_N$   must coincide with $\rho\vert_I$, which implies item (5).
\end{proof}

	\begin{proposition}
		Given $n\gg n_0>1$, for the word of length seven
		$$  w_{n,n_0} = (-v_n, -v_{n_0}, -v_{n_0}, v_{n_0}, -v_{n_0}, -v_{n_0}, -v_n ),    $$
		$A^7_t(w_{n,n_0})$ winds once around $\mathbb{P}^1$ as $t$ ranges in some interval
		\[  
		\left[ v^*_{n,n_0} - \frac{5}{ v_n},\, v^*_{n,n_0} + \frac{5}{v_n} \right] 
		\ \text{ with  } v^*_{n,n_0}  = v_{n_0} +   O\left(v_{n_0}^{-1}\right). 
		\]
	\end{proposition}

	\begin{proof}
		Let $\omega=(\cdots, -v_{n_0}, v_{n_0}, -v_{n_0}, \cdots)$ be the sequence from Proposition~\ref{prop hc asymp} below. We apply Proposition~\ref{prop:matchingsToTurns} with
		\[  \bld B(t) = (S_{ -v_{n_0}-t})^2  S_{ -v_{n}-t}   \ \text{ and } \   \bld A(t) = S_{ -v_{n}-t} \, (S_{ -v_{n_0}-t} ) ^2 S_{ v_{n_0}-t}  ,  \]		
	 	so that $(\bld A\,\bld B)(t) = A^7_t(w_{n, n_0})$.
		By Remark~\ref{rmk first block}, the lower bound on the winding speed of $\bld B(t)$ and $\bld A(t)^{-1}$ near $t=v_{n_0}$ follows from Proposition~\ref{prop lower bound winding speed} applied to the product of two Schr\"odinger matrices close to $S_{ -2v_{n_0}}$. That proposition gives, for $|t-v_{n_0}|\lesssim 1/v_{n_0}$, 
		\[
		c := \frac{2}{5 v_{n_0}^2} < \frac{2}{(2v_{n_0})^2} - O\!\left(\frac{1}{v_{n_0}^4}\right)
		\]
		as a lower bound on $|\text{winding speed}|$ for both $\bld A(t)^{-1}  \hat e_2$ and $\bld B(t) \, \hat e_1$.
		
		Within this range, a simple calculation yields
		\begin{align*}
			& 5 \, v_{n_0}^2 v_n   > 4 \, v_{n_0}^2\, v_n + O(v_{n_0}^2) \ge \norm{\bld B(t)} \\
			& \quad \ge \norm{\bld B(t)e_1} \ge (4 \, v_{n_0}^2-1)\, v_n - O(v_{n_0}^2) > 3 \, v_{n_0}^2 v_n,
		\end{align*}
		and
		\begin{align*}
			&5 \, v_{n_0}^2 v_n   > 4 \, v_{n_0}^2\, v_n + O(v_{n_0}^2) \ge \norm{\bld A(t)^{-1}}  \\
			& \quad \ge \norm{\bld A(t)^{-1}e_2} \ge  (4 \, v_{n_0}^2-1)\, v_n - O(v_{n_0}^2) > 3 \, v_{n_0}^2 v_n.
		\end{align*}  
		For $|t-v_{n_0}|\lesssim 1/v_{n_0}$, the leftmost matrix in $\bld A(t)^{-1}$ satisfies
		\begin{equation}\label{eq 90 graus}
		 	A_t(\omega)^{-1} = (S_{ v_{n_0}-t})^{-1} = \begin{bmatrix}  0 &  1 \\ -1 & O(1/v_{n_0}) \end{bmatrix} , 
		\end{equation}
		which is essentially a $90^{\circ}$ rotation with norm close to $1$.
		
		Set $\gamma = \log(3\, v_{n_0}^2 v_n)$ and $\gamma' = \log(5 \, v_{n_0}^2 v_n)$. To apply Proposition~\ref{prop:matchingsToTurns}, we show that for some $t_0 = v_{n_0} + O(1/v_{n_0})$, 
		\[
		\bld A(t_0)\,\bld B(t_0)\, \hat e_1 = \hat e_2,
		\]
		so that $(\bld B(t_0),\bld A(t_0))$ has a $(\gamma,\gamma')$-matching (Definition~\ref{def:matching}). The conclusion then follows because $5 c^{-1} e^{-\gamma} =  25/(6\,v_n)<5/v_n$.
		
		By Proposition~\ref{prop hc asymp}, there exists $t_* = v_{n_0} + 1/v_{n_0} + O(1/v_{n_0}^3)$ such that
		\[
		(S_{v_{n_0}-t_*})\,\hat w_+(v_{n_0}+ t_*) = \hat w_-(v_{n_0} + t_*), 
		\]
		where 
		\[  w_+(v_{n_0}+t_*) = (1, \lambda_+(v_{n_0}+t_*)^{-1}) \ \text{ and } \ w_-(v_{n_0}+t_*) = (\lambda_-(v_{n_0}+t_*), 1)    \]
		 are the unstable and stable eigenvectors of $(S_{-v_{n_0}-t_*})^2$.  
		 
		 Since $\lambda_-(v_{n_0}+t_*) = \lambda_+(v_{n_0}+t_*)^{-1} \sim -1/(2v_{n_0})$, we have
		\[
		d( S_{-v_n-t_*}\, \hat e_1, \hat w_+(v_{n_0}+t_*)) < \frac{1}{v_{n_0}}, \qquad
		d( (S_{-v_n-t_*})^{-1} \hat e_2, \hat w_-(v_{n_0}+t_*)) < \frac{1}{v_{n_0}} .
		\]
		Applying $(S_{-v_{n_0}-t_*})^2$ and $(S_{-v_{n_0}-t_*})^{-2}$  respectively to both projective points in the previous estimates,  we get that
		\begin{align*}
			d\bigl( \bld B(t_*)\,\hat e_1,\; \hat w_+(v_{n_0}+t_*)\bigr) &< \frac{1}{4v_{n_0}^3},\\
			d\bigl(\bld A(t_*)^{-1}\,\hat e_2,\; \hat w_-(v_{n_0}+t_*)\bigr) &< \frac{1}{4 v_{n_0}^3}.
		\end{align*}   
		Therefore
		\[
		d\bigl( \bld A(t_*)^{-1}\hat e_2,\; \bld B(t_*)\hat e_1\bigr)
		< \frac{1}{4v_{n_0}^3} + \frac{1}{4v_{n_0}^3} = \frac{1}{2v_{n_0}^3}.
		\]
		
		Both $\bld B(t)\, \hat e_1$ and $\bld A(t)^{-1} \hat e_2$ wind in opposite directions with absolute minimum speed $c>0$, so their relative speed is at least $2c$. Hence the equation
		\[
		\bld A(t)^{-1}\, \hat e_2 = \bld B(t) \, \hat e_1
		\]
		has a root $t_0 \approx t_*$ satisfying
		\[
		|t_0 - t_*| \le \frac{1/(2v_{n_0}^3)}{2c} = \frac{5v_{n_0}^2}{8 v_{n_0}^3} = \frac{5}{8 v_{n_0}}.
		\]
		Thus $t_0 = t_* + O(1/v_{n_0}) = v_{n_0} + O(1/v_{n_0})$, completing the proof.
	\end{proof}

Recall that given an interval  $J=[a,b]$ and a non-decreasing function $\rho(x)$ 
we write  $d\rho(J) =\rho(b)-\rho(a)$.

\begin{proof}[Proof of Theorem~\ref{main}] 
Fix $n \gg n_0 >1$ and consider the cylinder set $\Cscr_{n,n_0} \subset \R^\Z$ determined by the word 
\[
w_{n,n_0} = (-v_n, -v_{n_0}, -v_{n_0}, \, v_{n_0}, -v_{n_0}, -v_{n_0}, -v_n ). 
\]
Let $\nu, \nu_N \in \Prob(\R)$, $\nu=\lim_{N\to \infty} \nu_N$, where $\nu_N$ are the measures introduced in \eqref{eq nuN}.
By construction
\[
\nu_N^\Z(\Cscr_{n,n_0})  = \Bigl(\sum_{j=1}^N p_j\Bigr)^{-7} \frac{p_{n_0}^5\, p_n^2}{128},
\qquad  
\nu^\Z(\Cscr_n) = \frac{p_{n_0}^5\, p_n^2}{128}  
.
\]

Let $L \in \N$ be large and let $\omega \in \R^\Z$ be $\nu_N^\Z$-typical, in the sense of the Birkhoff Ergodic Theorem applied to the shift $\sigma^7 : \R^\Z \to \R^\Z$ and the indicator function of $\Cscr_{n,n_0}$.  
Define
\[
\Sigma_L := \{\, j \in \{0,1,\ldots,L-1\} : \sigma^{7j}(\omega) \in \Cscr_{n,n_0} \,\}. 
\]
Each index $j \in \Sigma_L$ corresponds to an occurrence of $w_{n,n_0}$-matching at some parameter $t \in J_{n,n_0}$, where
\[
J_{n,n_0} := \Bigl[v_{n,n_0}^\ast - \tfrac{5}{v_n},\; v_{n,n_0}^\ast + \tfrac{5}{ v_n}\Bigr].
\]
By Propositions~\ref{prop:matchingsToTurns} and~\cite[Proposition 3.18]{BCDFK}  this yields
\[
d\rho_N(J_{n,n_0}) \;\geq\; \lim_{L \to \infty} \frac{\#\Sigma_L}{7L}
= \Bigl(\sum_{j=1}^N p_j\Bigr)^{-7} \frac{p_{n_0}^5\, p_n^2}{896}. 
\]
Observe that, for any fixed $\hat v\in\Proj$, the cardinality $\#\Sigma_L$ provides a lower bound for the number of full turns around $\Proj$ made by the curve $A_t^{7L}(\omega)\hat v$ as $t$ ranges over $J_{n,n_0}$.

Now consider the modulus of continuity
\[
\omega(r) := \frac{1}{\varphi(\log(1/r))^2}, 
\]
which is at least $2$-log H\"older.  
Since $\omega(|J_{n,n_0}|) = \varphi(\log(1/|J_{n,n_0}|))^{-2}$, for $N \gg n$ we obtain
\[
\frac{d\rho_N(J_{n,n_0})}{\omega(|J_{n,n_0}|)} 
\;\geq\; 
\Bigl(\sum_{j=1}^N p_j\Bigr)^{-7} \frac{p_{n_0}^5\, p_n^2}{896} \, \varphi\!\left(\log \tfrac{1}{|J_{n,n_0}|}\right)^2. 
\]

By Item~(5) of Proposition~\ref{prop existence of ids}, we may pass to the limit as $N \to \infty$.  
Since $|J_{n,n_0}| = C/v_n$ for some $C > 0$ and $p_n \, \varphi(\log v_n) \to +\infty$ as $n \to \infty$, we deduce

\begin{align*}
	\frac{ d\rho(J_{n,n_0})}{\omega(|J_{n,n_0}|)} 
	&\geq \frac{p_{n_0}^5}{896} \, p_n^2\, \varphi\!\bigl(\log v_n - \log C \bigr)^2 \\
	&\gtrsim \bigl( p_{n-1}\, \varphi(\log v_{n-1}) \bigr)^2 \;\longrightarrow\; +\infty
	\qquad (n \to \infty).
\end{align*}

Thus, the function  $t \mapsto \rho(t)$ cannot have the local modulus of continuity $\omega$ in the interval $J_{n,n_0}\subset I_{n_0}$.
\end{proof}

\begin{proof}[Proof of Corollary \ref{maincoro}]
By the Thouless formula,  item (6) of Proposition~\ref{prop existence of ids}, \ $L_1(\mu_t)$ is the 
Hilbert transform of $\rho(t)$. By Theorem~\ref{main},  $\rho(t)$ does not have the modulus of continuity $\omega=\Tscr_2(\varphi)$.

\noindent
(a)\;  Notice that ~\cite[Lemma 10.3]{GS01} 
applies to moduli of continuity satisfying  equations (10.9) and (10.10) of 
\cite{GS01}, which includes H\"older and weak-H\"older. Thus, when $\varphi(r)$ is exponential or subexponential, by ~\cite[Lemma 10.3]{GS01} $L_1(\mu_t)$ can not have the modulus of continuity $\omega$. 

\noindent
(b)\;  If  $\varphi(r)=r^p$, for some $p>1$,  then by Theorem~\ref{main}  $\rho(t)$ is not $(2p)$-log-H\"older continuous. Hence by ~\cite[Corollary 2.3]{ALSZ21}, $L_1(\mu_t)$ can not be $(2p-1)$-log-H\"older continuous.  
\end{proof}

\begin{proof}[Proof of Corollary \ref{corollary}]
Apply Theorem~\ref{main} with 
\[
\varphi(r) = e^{r^{2/3}}, 
\qquad 
\psi(r) = e^{r^{1/3}}, 
\qquad 
p_n = \frac{6}{\pi^2 n^2}, 
\qquad 
v_n = \exp\!\bigl( (3 \log n)^{3/2} \bigr).
\]

We first verify the hypotheses of the theorem.  
Clearly $\sum_{n \geq 1} p_n = 1$, and
\[
\limsup_{n \to \infty} \frac{p_{n-1}}{p_n} 
= \lim_{n \to \infty} \frac{n^2}{(n-1)^2} 
= 1 ,
\]
so condition (2) holds.

A numerical calculation gives $v_2-v_1>4$. 
Since $\lim_{n \to \infty} v_n = +\infty$, it remains to show that 
\[
\lim_{n \to \infty} (v_n - v_{n-1}) = +\infty .
\]
Define $f(x) := \exp\!\bigl( (3 \log x)^{3/2} \bigr)$ for $x \geq 2$.  
Then $f$ is $C^1$ and strictly increasing.  
By the Mean Value Theorem, for each $n \geq 3$ there exists $\xi_n \in (n-1,n)$ such that
\[
v_n - v_{n-1} = f(n) - f(n-1) = f'(\xi_n).
\]
We compute
\[
f'(x) 
= f(x)\, \frac{d}{dx}\bigl( (3 \log x)^{3/2} \bigr) 
= \frac{9 \sqrt{3}}{2}\, \frac{\sqrt{\log x}}{x}\, 
\exp\!\bigl( (3 \log x)^{3/2} \bigr).
\]
As $x \to \infty$, $f'(x) \to +\infty$.  
Since $\xi_n \to \infty$, we conclude that
\[
v_n - v_{n-1} = f'(\xi_n) \;\xrightarrow[n \to \infty]{}\; +\infty,
\]
and condition (3) follows.

Given $t_0\in \R$, there exists $n_0=n_0(t_0)\geq 1$ such that for all $|t-t_0|<1$ and $n\geq n_0$,
\begin{align*}
\log \norm{A_{\pm v_n, t}} &\leq \log(v_n +|t|+1) = \log \left(1+\frac{|t|+1}{v_n} \right) + \log v_n \\ &\leq  \frac{|t|+1}{v_n} +\log v_n \leq  1+ \log v_n
\end{align*} 
and so for $n$ large enough, say $n\geq N_0>n_0$,
\begin{align*}
	\psi\left( \log \norm{A_{\pm v_n, t}}\right) &\leq \psi\left(  1+ \log v_n \right)
	= e^{(1+\log v_n)^{1/3}} \\
	& \leq e^{2\, (\log v_n)^{1/3}} \leq e^{2\, \sqrt{ 3\, \log n }} \leq e^{\frac{1}{2}\, \log n } =n^{\frac{1}{2}}
\end{align*} 
which implies that for all $t\in\R$ such that  $|t-t_0|<1$, 
\begin{align*}
\int \psi(\log \norm{g})\, d\mu_t(g) &\leq  \sum_{n \geq 1} p_n \psi(1 + \log v_n  ) \\
& \leq C + \frac{6}{\pi^2} \sum_{n \geq N_0} \frac{1}{n^{3/2}} < \infty .
\end{align*} 
Constant $C$ is the sum of  the first $N_0-1$ terms.
This proves that condition (4) holds.  
  
We compute
\[
\lim_{n \to \infty} p_n \varphi(\log v_n) 
= \lim_{n \to \infty} \frac{6}{\pi^2 n^2} \, e^{(\log v_n)^{2/3}}.
\]
Since $(\log v_n)^{2/3} = 3 \log n$, this becomes
\[
\lim_{n \to \infty} \frac{6}{\pi^2 n^2} \, e^{3 \log n}
= \lim_{n \to \infty} \frac{6}{\pi^2} \, n = +\infty .
\]
Thus condition (5) is satisfied.

We have therefore verified conditions (2)--(5) of Theorem~\ref{main}.  
It follows that the Lyapunov exponent function 
\[
\mathbb{R} \ni t \mapsto L_1(\mu_t)
\]
cannot have the local modulus of continuity
\[
\omega(r) = \bigl(\varphi(\log(1/r))\bigr)^{-2} 
= e^{-2 (\log(1/r))^{2/3}}. 
\]
In particular, $L_1(\mu_t)$ is not $(2,2/3)$-weak-H\"older continuous, and hence is not $\alpha$-Hölder continuous for any $\alpha > 0$.
\end{proof}

\begin{remark}
The previous argument provides insight into the fractal character of $\rho(t)$ (the IDS). In each component $I_{n_0}$ of the spectrum
 $\Sigma$, the regularity of $\rho(t)$ is determined by higher energy values of the potential far away from $I_{n_0}$. The proof uses a particular sequence of matchings $w_{n, n_0}$ attained at specific energy values $t=v^\ast_{n, n_0}\in I_{n_0}$, but the same irregularity of $\rho(t)$ should also arise from different sequences of matchings, leading to different energy values in $I_{n_0}$ where regularity breaks down. This naturally raises Question~\ref{second question}.
\end{remark}

\bigskip

\section{Appendix: Some linear algebra facts}

This appendix recalls definitions from~\cite{BCDFK} for the reader's convenience and extends a particular result from that work, which is used in the proof of Theorem 2.1.

\subsection*{Winding property}

Consider a smooth family  
$\bld A:I\to \SL_2(\R)$ defined on some interval $I\subset \R$.

Given any two vectors $v_1, v_2\in\R^2$, recall  that
$$v_1\wedge v_2 := \det(v_1, v_2)=\norm{v_1} \norm{v_2} \sin \measuredangle(v_1, v_2). $$

\begin{definition}
	\label{winding def}
	We say that a curve $\bld A$ is positively (respectively negatively)  winding if for all $t\in I$, the quadratic form
	$Q_{\bld A, t}:\R^2\to\R$,
	$$Q_{\bld A, t}(v):= (\bld A(t)  \,v)\wedge (\dot {\bld A}(t)\, v) =    v \wedge ( {\bld A} (t)^{-1} \, \dot {\bld A}(t) \, v),$$
	is  positive (respectively negative) semidefinite, with one eigenvalue bounded away from $0$, where $\dot {\bld A}(t):=  \frac{d \bld A}{dt}(t)$. Moreover we say that $\bld A$ is strictly positively (respectively negatively)  winding if for all $t\in I$, the quadratic form  $Q_{\bld A, t}$ is positive  (respectively negative) definite for all $t\in I$  with eigenvalues uniformly bounded away from $0$.
\end{definition}

The name and relevance of this quadratic form emerge from the following formula (see~\cite[Proposition 3.1]{BCDFK}):

\[ \frac{d}{dt} \frac{\bld A(t)\, v}{\norm{\bld A (t)\, v}} = \left[\, \frac{ Q_{\bld A, t}(v) }{ \norm{\bld A (t)\, v}^2 }  \,\right]\, \bld z ( \bld A (t)\, v ) ,  \]
where $\bld z(w)$ denotes the unique
unit vector which makes $\{\frac{w}{\norm{w}}, \bld z(w) \}$ a positively oriented  orthonormal basis,  for any non-zero vector $w \in \R^2$.

This formula shows that the positive or negative winding character of $\bld A$ translates to all curves
$I \ni t \mapsto \frac{\bld A(t)\, v}{\norm{\bld A (t)\, v}}$ (with values in the unit circle or in projective space) having speeds that do not change sign.

\begin{remark}
	\label{rmk first block}
	From \cite[Proposition 3.4]{BCDFK} it follows that, given smooth curves $\bld A, \bld B : I \to \SL_2(\R)$, if $\bld B$ is positively winding while $\bld A$ is strictly positively winding with a uniform lower bound $c>0$ on the speeds of the curves $I \ni t \mapsto \frac{\bld A(t)\, v}{\norm{\bld A (t)\, v}}$ (i.e., for all $t\in I$ and every unit vector $v\in \R^2$), then the composition curve $(\bld A\, \bld B)(t) = \bld A(t)\, \bld B(t)$ is also strictly positively winding, with the same uniform lower bound $c>0$ on the speeds of its induced curves $I \ni t \mapsto \frac{ (\bld A \bld B)(t)\, v}{\norm{\bld A \bld B (t)\, v}}$.
\end{remark}

We use below the Schr\"odinger matrix notation introduced in~\eqref{Sch mat}.
\begin{example}
	The Schr\"odinger family $\bld A(t)=S_{v-t}$ is positively (but not strictly positive) winding because
	$$ Q_{\bld A, t}(x,y) = x^2\geq 0 . $$	
\end{example}

\begin{example}
	The Schr\"odinger product $\bld B(t)=S_{v_1-t}\, S_{ v_2-t } $
	is strictly positively  winding with
	$$ Q_{\bld B, t}(x,y) = \begin{bmatrix}
		x & y 
	\end{bmatrix}\,
	\begin{bmatrix}
		2+ 2 (t-v_2)^2 & 2 (t-v_2)\\
		2 (t-v_2) & 2 \\
	\end{bmatrix}\,  \begin{bmatrix}
		x \\ y 
	\end{bmatrix} . $$	
	Notice that the matrix associated with this quadratic form has determinant $=4$ and trace  $=4+2 (t-v_2)^2\geq 4$, which implies the positive definiteness of the matrix.
\end{example}

\begin{proposition} 
	\label{prop lower bound winding speed}
	The Schr\"odinger family 
	\[ {\bld B}_v(t)= ( S_{-v-t})^2  =\begin{bmatrix}
		(v+t)^2-1 & v+t \\ -(v+t) & -1
	\end{bmatrix}  \]
	admits the following lower bound on its minimum winding speed
	\[  \min_{\norm{(x,y)}=1} \frac{ Q_{{\bld B}_v, t}(x,y) }{ \norm{{\bld B}_v (t)\, (x,y)}^2 }   \geq \frac{2}{(v+t)^2}-O\left( \frac{1}{(v+t)^4}\right) .  \]
\end{proposition}

\begin{proof}
	For $v>0$,  the quadratic form
	\[ Q_v(x,y)= \left( (S_{-v})^2 \begin{bmatrix}
		x\\y
	\end{bmatrix}\right) \wedge  \left(\frac{d}{dv}(S_{-v})^2 \begin{bmatrix}
		x\\y
	\end{bmatrix}\right)\]
	is defined by the positive definite symmetric matrix
	\[  M(v)=\begin{bmatrix}
		2+ 2 v^2 & 2 v\\
		2 v & 2 \\
	\end{bmatrix} . \]
	Similarly, the quadratic form
	\[ P_v(x,y)= \left\| (S_{-v})^2 \begin{bmatrix}
		x\\y
	\end{bmatrix} \right\|^2  \]
	is defined by the positive definite symmetric matrix
	\[  N(v)=\begin{bmatrix}
		v^4-v^2+1 & v^3 \\
		v^3 & v^2+1 \\
	\end{bmatrix} .\]
	An easy calculation shows that the minimum value
	\[ \min_{\norm{(x,y)}=1} \frac{Q_v(x,y)}{P_v(x,y)}  \]
	corresponds to the minimum eigenvalue solution of the generalized eigenvalue problem
	\[   M_v \begin{bmatrix}
		x\\y
	\end{bmatrix} =\lambda\, N_v \begin{bmatrix}
		x\\y
	\end{bmatrix}\quad \Leftrightarrow\quad (N_v)^{-1} M_v \begin{bmatrix}
		x\\y
	\end{bmatrix} =\lambda\,  \begin{bmatrix}
		x\\y
	\end{bmatrix} \]
	Finally, since 
	\[ (N_v)^{-1} M_v = \begin{bmatrix}
		4 v^2+2 & 2 v \\
		2 v-4 v^3 & 2-2 v^2 \\
	\end{bmatrix} \]
	has minimum eigenvalue
	\[ \lambda_{min}= 2 + v^2-v \,\sqrt{v^2+4} =  \frac{2}{v^2}-O\left( \frac{1}{v^4}\right)\]
	substituting $v$ by $v+t$ we get the stated lower bound.
\end{proof}

\begin{proposition}
	\label{prop hc asymp}
	Consider 
	$\omega=(\omega_j)_{j\in \Z}\in\R^\Z$   with    $\omega_0=v>0$   and $\omega_j= 
	-v$  for $j\neq 0$.  
	Then the Schr\"odinger operator $H_\omega:\ell^2(\Z)\to \ell^2(\Z)$  with potential $\omega$  has an eigenvalue at  some\ $t_\ast=v+\frac{1}{v} +O\left(\frac{1}{v^3}\right)$. 
	
	In other words, the Schr\"odinger cocycle  
	$A_t:\R^\Z\to \SL_2(\R)$,      $A_t(\omega)=S_{ \omega_0-t }$,  has a  heteroclinic connection, i.e., $E^u_{A_t}(\omega)=E^s_{A_t}(\omega)$,  
	for this particular sequence $\omega$ at some \ $t_\ast=v+\frac{1}{v} +O\left(\frac{1}{v^3}\right)$.
\end{proposition}

\begin{proof}
	The eigenvalues of   $S_{-v}=\begin{bmatrix}
		-v & -1\\
		1 & 0 \\
	\end{bmatrix}$ \ ($v>0$) \
	are 
	\[\lambda_\pm(v) = \frac{1}{2} \left(-v\pm \sqrt{v^2-4} \right) \ \text{ with } \ \lambda_-(v)<-1<\lambda_+(v)<0 \] 
	and the corresponding eigenvectors are
	\[ w_\pm (v)= \left( \lambda_\pm(v), 1 \right) . \]
	Hence the required parameter $t=t_\ast$ must be such that  $w_-(v+t_\ast)$ is collinear with $(S_{v-t_\ast})\, w_+(v+t_\ast)$.
	Take 
	\[ t_\ast = -v + 2 \sqrt{v^2+1} =v+\frac{1}{v} +O\left( \frac{1}{v^3}\right)  .\] 
	Then $v+t_\ast= 2 \sqrt{v^2+1}$ and 
	$-v+t_\ast= -2 v + 2 \sqrt{v^2+1}$. Hence
	\begin{align*} 
		\lambda_\pm(v+t_\ast) &= \mp v - \sqrt{v^2+1}, \\
		w_\pm(v+t_\ast) &= \left( \mp v - \sqrt{v^2+1} , \, 1 \right)  ,
	\end{align*}
	while
	\begin{align*}  
		 S_{v-t_\ast}  \begin{bmatrix}
			\lambda_+(v+t_\ast)\\ 1
		\end{bmatrix}     &=\begin{bmatrix}
			2 v - 2 \sqrt{v^2+1} & -1\\ 1 & 0
		\end{bmatrix} \begin{bmatrix}
			- v - \sqrt{v^2+1}\\ 1
		\end{bmatrix} \\
		&=\begin{bmatrix}
			1\\- v - \sqrt{v^2+1} 
		\end{bmatrix} =
		\begin{bmatrix}
			1\\\lambda_+(v+t_\ast) 
		\end{bmatrix} \\
		&= \lambda_+(v+t_\ast) \, \begin{bmatrix}
			\lambda_-(v+t_\ast) \\ 1
		\end{bmatrix}
	\end{align*}
	Therefore  the stated result holds for  $t= t_\ast$.
\end{proof}

\bigskip

\subsection*{ Singular directions and projective action }

Given $A\in\SL_2(\R)$
denote by $\{\medir(A),\, \ledir(A)\}$  an orthonormal set of singular directions of $A$,  defined by the relations 
\begin{align*}
	(A^\ast\,A)\, \medir(A)=\norm{A}^2\, \medir(A)
	\quad
	\text{ and }
	\quad
	(A^\ast\,A) \,\ledir(A)=\norm{A}^{-2}\,\ledir(A) ,  
\end{align*}
where  $A^\ast$ denotes the transpose of $A$.
Define also
\begin{align*}
	\medir^\ast(A) := \medir(A^\ast)
	\quad
	\text{and}
	\quad
	\ledir^\ast(A) := \ledir(A^\ast),
\end{align*}
so that
\begin{align*}
	A\, \medir(A)=\norm{A}\,\medir^\ast(A)
	\quad
	\text{and}
	\quad
	A\,\ledir(A)= \norm{A}^{-1}\,\,\ledir^\ast(A).
\end{align*}

For each pair of vectors $v,\, w\in \R^2$ we write 
\begin{align*}
	d(\hat v,\, \hat w) := \frac{|v\wedge w|}{\norm{v}\, \norm{w}} = \sin \measuredangle(v, w),
\end{align*}
for the usual distance in $\mathbb{P}^1$.
\begin{lemma}[\cite{BCDFK}, Lemma 6.1] \label{L1}
	Given $A\in \SL_2(\R)$, for every $\hat v\in \mathbb{P}^1$,
	\begin{align*}
		\frac{1}{\norm{A}}\sqrt{
			\norm{A\, v}^2 - \norm{A}^{-2}
		}
		\leq d(\hat v,\, \hat \ledir(A))
		\leq \frac{\norm{A\, v}}{\norm{A}} \ .
	\end{align*}
\end{lemma}

\begin{lemma}[\cite{BCDFK}, Lemma 6.2]\label{L2} 
	Given $A\in \SL_2(\R)$, for every $\hat v\in \Pp^1$,
	\begin{align*}
		d(A\hat v,\, \hat \medir^*(A))  \leq \frac{1}{d(\hat v,\,  \hat \ledir(A))\, \norm{A}^2 }  .
	\end{align*}
\end{lemma}

\subsection*{ Matchings }

The following definition extends the notion of matching introduced in~\cite[Definition 5.1]{BCDFK}.
\begin{definition}
	\label{def:matching}
	Given $\gamma'>\gamma> 0$, we say that a pair of matrices $(B,\, A)\in \SL_2(\R)^2$ has a $(\gamma, \gamma')$-matching if there exists a pair of points  $\hat e_1\, \hat e_2 \in \mathbb{P}^1$ such that
	\begin{enumerate}
		\item $A\, B\, \hat e_1 = \hat e_2$;
		
		\item $e^{\gamma}
		\leq \norm{B\, e_1}
		\leq \norm{B} 
		\leq 2\, e^{\gamma'}$;
		
		\item $e^{\gamma}
		\leq \norm{A^{-1}\, e_2}
		\leq \norm{A} 
		\leq 2\, e^{\gamma'}$.
	\end{enumerate}
	In this case we say that $(B,\, A)$ connects $\hat e_1,\, \hat e_2 \in \mathbb{P}^1$.
\end{definition}

\begin{proposition}\label{prop:matchingsToAlmostTurn}
	There exists $\gamma_0>0$ such that for every $\gamma' >\gamma \geq \gamma_0$, if $(B, A)\in \SL_2(\R)^2$ is a $(\gamma, \gamma')$-matching and $\hat v,\, \hat w\in \mathbb{P}^1$ satisfy
	\begin{align*}
		d(\hat v,\, \hat \ledir(B)) \geq e^{-\gamma}
		\quad
		\text{and}
		\quad
		d(\hat w,\, \hat \medir^*(A)) \geq e^{-\gamma},
	\end{align*}
	then
	\begin{align*}
		d(B\, \hat v,\, A^{-1}\, \hat w) \leq 3\, e^{-\gamma}.
	\end{align*}
\end{proposition}
\begin{proof}
	Assume $(B, A)$ connects $\hat e_1$, $\hat e_2\in\Pp^1$. By Lemma \ref{L2} and the matching assumption,
	\begin{align*}
		d(B\, \hat v,\, \hat \medir^*(B)) \leq \frac{1}{d(\hat v,\, \hat \ledir(B))\norm{B}^2} \leq e^{-\gamma}
	\end{align*}
	and
	\begin{align*}
		d(A^{-1}\, \hat w,\, \hat \ledir(A)) \leq \frac{1}{d(\hat w,\, \hat \medir^*(A))\norm{A}^2} \leq e^{-\gamma}.
	\end{align*}
	Using both lemmas \ref{L1} and~\ref{L2},
	\begin{align*}
		d(\hat \medir^*(B),\, B\, \hat e_1)
		&\leq\frac{1}{d(\hat e_1,\, \hat \ledir(B))\norm{B}^2}
		\leq \frac{\norm{B}}{\norm{B}^2\, \sqrt{\norm{B\, e_1}^2 - \norm{B}^{-2}}}\\
		&\leq e^{-\gamma}\frac{1}{\sqrt{e^{2\gamma} - e^{-2\gamma'}}} \leq 2\, e^{-2\gamma},
	\end{align*}
	and similarly,
	\begin{align*}
		d(\ledir(A),\, A^{-1}\, \hat e_2) \leq \, 2\, e^{-2\gamma},
	\end{align*}
	for every $\gamma\geq \gamma_0$, for some $\gamma_0$ sufficiently large.
	
	Thus, by the triangle inequality and the fact that $B\, \hat e_1 = A^{-1}\, \hat e_2$,
	\begin{align*}
		d(B\, \hat v,\, A^{-1}\, \hat w)
		&\leq d(B\, \hat v,\, \hat \medir^*(B)) + d(\hat \medir^*(B),\, B\, \hat e_1)\\
		&+ d(A^{-1}\, \hat e_2,\, \hat \ledir(A)) + d(\hat \ledir(A),\, A^{-1}\, \hat w)\\
		&\leq 2\,(e^{-\gamma} + 2e^{-2\gamma}) \leq 3\, e^{-\gamma},
	\end{align*}
	for every $\gamma\geq \gamma_0$.
\end{proof}

\bigskip

\begin{definition} 
\label{def length}
Given  $I\ni t\mapsto A_t\in \SL_2(\R)$ smooth and a unit vector $v\in\Su^1$,  the length (oriented angle) of $A_t\, v/\norm{A_t\, v}\in\Su^1$ as $t$ ranges in $I$ is denoted by $\ell_I(A_t \, v)$.
	This also agrees with the length of the projective curve $I\ni t\mapsto A_t\,\hat v$.
\end{definition}

\begin{definition}
We say that a smooth curve $I\ni t\mapsto A_t\in\SL_2(\R)$ winds $n$ times around $\Pp^1$ if
	for all $\hat v\in\Pp^1$, $\ell_I(A_t\, \hat v)\geq n\,\pi$.
\end{definition}

The following proposition generalizes~\cite[Proposition 5.6]{BCDFK}.

\begin{proposition} 
	\label{prop:matchingsToTurns}
	There exists $\gamma_0>0$ such that for any $c>0$ and for every pair of smooth, positively strictly  winding families of matrices $\{\bld A_t\}_{t\in I}$ and $\{\bld B_t\}_{t\in I}$ with winding speed bounded from below by $c$,  if
	$(\bld B_{t_0},\, \bld A_{t_0})\in \SL_2(\R)^2$ has a $(\gamma, \gamma')$-matching where $\gamma'>\gamma\geq \gamma_0$ then
	$\bld A_t\, \bld B_t$ winds once around $\mathbb{P}^1$ as $t$ ranges in $[t_0 - 5\, c^{-1}\, e^{-\gamma},\, t_0 + 5\, c^{-1}\, e^{-\gamma}]$.
\end{proposition}
\begin{proof}
	Consider  directions $\hat e_1,\, \hat e_2 \in \mathbb{P}^1$ which are connected by 	$(\bld B_{t_0},\, \bld A_{t_0})$. Set
	\begin{align*}
		J_0 := [t_0 - 3 \, c^{-1}\, e^{-\gamma},\, t_0 + 3 \,c^{-1}\, e^{-\gamma}]
	\end{align*}
	
	For any pair of vectors $\hat v,\, \hat w\in \mathbb{P}^1$ satisfying
	\begin{align*}
		d(\hat v,\, \hat \ledir(\bld B_{t_0})  ) \geq e^{-\gamma}
		\quad
		\text{and}
		\quad
		d(\hat w,\,  \hat \medir^\ast(\bld A_{t_0} ) ) \geq e^{-\gamma},
	\end{align*}
	we can use Proposition \ref{prop:matchingsToAlmostTurn} to conclude that 
	\[ d(\bld B_{t_0}\, \hat v,\, \bld A_{t_0}^{-1}\, \hat w) \leq 3\, e^{-\gamma} . \] 
	If we define $f_+,\, f_-:J_0\to \mathbb{P}^1$ by,
	\begin{align*}
		f_+(t) := \bld B_t\, \hat v
		\quad
		\text{and}
		\quad
		f_-(t) := \bld A_t^{-1}\, \hat w,
	\end{align*}
	then by the  strict winding hypothesis, see~ Definition~\ref{winding def},  we have that $f_-'(t) < -c < 0 < c < f_+'(t)$, for every $t\in  I$. 
	
	Since $d(f_+(t_0),\, f_-(t_0)) \leq 3\, e^{-\gamma}$, there exists $t_*\in J_0$ such that
	\begin{align*}
		\bld B_{t_\ast}\, \hat v
		= f_+(t_\ast) = f_-(t_\ast)
		= \bld A_{t_\ast}^{-1}\, \hat w.
	\end{align*}
	Using Corollary~\cite[Corollary 3.16]{BCDFK}, this implies	
	that  given $\hat v\in \Pp^1$ with 	$d(\hat v,\, \hat \ledir(\bld B_{t_0})  ) \geq e^{-\gamma}$,	 
	\begin{align*}
		\ell_{J_0}( \bld A_t\, \bld B_t \, \hat v) \geq \pi - 2\, e^{-\gamma}.
	\end{align*}
	Writing $J_1 := J_-\cup J_0 \cup J_+$, where $J_- := [t_0 - 4\, c^{-1}\, e^{-\gamma},\, t_0 - 3 c^{-1}\, e^{-\gamma}]$ and 
	$J^+ := [ t_0 + 3 c^{-1}\, e^{-\gamma},\, t_0 + 4\, c^{-1}\, e^{-\gamma}]$ and using again the strict positive winding, we have
	\begin{align*}
		\ell_{J_1}(\bld A_t\, \bld B_t\, \hat v)
		&\geq \ell_{J_+}(\bld A_t\, \bld B_t\, \hat v) + \ell_{J_0}(\bld A_t\, \bld B_t\, \hat v) + \ell_{J_-}(\bld A_t\, \bld B_t\, \hat v) \\
		&\geq  2 \, e^{-\gamma} + \pi - 2 \,e^{-\gamma} = \pi.
	\end{align*}
	In other words, for every $\hat v\notin \text{B}(\hat \ledir(\bld B_{t_0}), e^{-\gamma})$, the curve $t\in J_1\mapsto \bld A_t\, \bld B_t \, \hat v$ gives one turn around $\mathbb{P}^1$. Thus, by~\cite[Corollary 3.16]{BCDFK} for every $\hat w\in \mathbb{P}^1$ and $\hat v\notin \text{B}(\hat \ledir(\bld B_{t_0}), e^{-\gamma})$, the equation $\bld B_t^{-1}\, \bld A_t^{-1}\,   \hat w = \hat v$ has a solution for some $t\in J_1$ which implies again by~\cite[Corollary 3.16]{BCDFK} that
	\begin{align*}
		\ell_{J_1}(\bld B_t^{-1}\, \bld A_t^{-1}\, \hat w) \geq \pi - 2 e^{-\gamma}.
	\end{align*}
	Taking $J_* := [t_0 - 5 \,c^{-1}\, e^{-\gamma},\, t_0 + 5 \, c^{-1}\, e^{-\gamma}]$ and using the same argument as above we conclude that for every $\hat w\in \mathbb{P}^1$,
	\begin{align*}
		\ell_{J_*}(\bld B_t^{-1}\, \bld A_t^{-1}\,  \hat w) \geq \pi.
	\end{align*}
	Therefore, for every $\hat v,\, \hat w\in \mathbb{P}^1$, the equation
	\begin{align*}
		\bld B_t^{-1}\, \bld A_t^{-1}\, \hat w = \hat v
		\quad \Leftrightarrow \quad
		\quad
		\bld A_t\, \bld B_t\, \hat v = \hat w,
	\end{align*}
	has a solution for some $t\in J_*$, and an application of~\cite[Corollary 3.17]{BCDFK} concludes the proof.
\end{proof}

\bigskip



%
%
%
%
%
%
%
%
%


\providecommand{\bysame}{\leavevmode\hbox to3em{\hrulefill}\thinspace}
\providecommand{\MR}{\relax\ifhmode\unskip\space\fi MR }
\providecommand{\MRhref}[2]{%
	\href{http://www.ams.org/mathscinet-getitem?mr=#1}{#2}
}
\providecommand{\href}[2]{#2}

\end{document}